# The general solution of the linear difference equation of degree-2 and the continued fraction made from this equation

Nikos Bagis
Department of Informatics
Aristotle University of Thessaloniki
Greece
Email: bagkis@hotmail.com

**Abstract**
In this article we give, for the fist time the solution of the general difference equation of 2-degree. We also give as application the expansion of a continued fraction into series, which was first proved, found in the past by the author.

**Introduction**

It is well known that the linearly difference sequence

$$x_{n+2} = b_n x_{n+1} + a_n x_n \quad :(1)$$

produces a continued fraction $K$ .i.e.

$$K = \cfrac{a_0}{b_0 + \cfrac{a_1}{b_1 + \cfrac{a_2}{b_2 + \dots}}} = \underset{n=1}{\overset{\infty}{K}}\left(\frac{a_n}{b_n}\right) \quad :(2)$$

It is also well known that the tail sequence for this fraction is:

$$t_n = \frac{a_n}{b_n + t_{n+1}}, \quad n \in \mathbb{N}, \quad t_n = -\frac{x_{n+1}}{x_n} \quad :(3)$$

and if $(x_n)$ is a "minimal" solution of (1) see [L,W] , then $t_0$=K.

**Theorems**

**Theorem 1.**
The general solution of (1) is :

$$\frac{x_n}{\prod_{i=1}^{n-2} b_i} = x_1 \cdot \left(1 + \sum_{i=1}^{n-2} g_i + \sum_{S_{1,n}(i,j)} g_i g_j + \sum_{S_{1,n}(i,j,k)} g_i g_j g_k + \dots\right) +$$

$$+ x_0 \cdot \frac{a_0}{b_0} \cdot \left(1 + \sum_{i=2}^{n-2} g_i + \sum_{S_{2,n}(i,j)} g_i g_j + \sum_{S_{2,n}(i,j,k)} g_i g_j g_k + \dots\right)$$

:(4)

Where:
$$S_{q,n}(i,j,k,...,l,m) =$$
$$= \{i,j,...,m \in N : q \le i < j < k < ... < l < m \le n-2, j-i \ge 2, k-j \ge 2, ..., m-l \ge 2\}$$

$q = 1, 2$, and $g_i = \dfrac{a_i}{b_{i-1}b_i}$, $i = 1, 2, 3, ...$

Knowing that $x_n$ is a zero limit sequence we take the limits in (4) and calculate the ratio $-x_1/x_0 = t_0 = K$. After some calculations we can see that:

**Theorem 2.**

$$\cfrac{a_0}{b_0 + \cfrac{a_1}{b_1 + \cfrac{a_2}{b_2 + ....}}} = \frac{a_0}{b_0} \cdot \frac{1 + \sum_{i=2}^{\infty} g_i + \sum_{i=2}^{\infty}\sum_{j=1}^{\infty} g_i g_{i+j+1} + \sum_{i=2}^{\infty}\sum_{j=1}^{\infty}\sum_{k=1}^{\infty} g_i g_{i+j+1} g_{i+j+k+2} + ...}{1 + \sum_{i=1}^{\infty} g_i + \sum_{i=1}^{\infty}\sum_{j=1}^{\infty} g_i g_{i+j+1} + \sum_{i=1}^{\infty}\sum_{j=1}^{\infty}\sum_{k=1}^{\infty} g_i g_{i+j+1} g_{i+j+k+2} + ...} \quad :(5)$$

Without the care of convergence.

(For another proof see [Ba]).

At this point we can see some applications of this continued fraction expansion.

## Applications

**Application 1**

For $a_m = z$ and $b_m = m+c$, $z \in C$ we find:

$$\cfrac{c}{c + \cfrac{z}{c+1 + \cfrac{z}{c+2 + \cfrac{z}{c+3 + \cfrac{z}{c+4 + ....}}}}} = \frac{{}_0F_1(c+1;z)}{{}_0F_1(c;z)}$$

Where ${}_0F_1$ the usual Hypergeometric function.

**Application 2**

For $a_m = q^m$ and $b_m = 1$ we find the Rogers-Ramanujan continued fraction expansion:

$$\cfrac{1}{1+\cfrac{qz}{1+\cfrac{q^2 z}{1+\cfrac{q^3 z}{1+...}}}} = \frac{\sum_{k=0}^{\infty} \frac{q^{k(k+1)} z^k}{(q)_k}}{\sum_{k=0}^{\infty} \frac{q^{k^2} z^k}{(q)_k}}$$

where:
$$(q)_k = (1-q)\cdot(1-q^2)\cdot...\cdot(1-q^k), \quad (q)_0 = 1$$

**Application 3**

For $a_m = z$ and $b_m = q^{(-1)^m \left(c + \sum_{k=0}^{m}(-1)^k k\right)}$, $c \in \mathbb{R}$, we find (after some elementary calculations)

$$\cfrac{1}{1+\cfrac{z}{q+\cfrac{z}{q+\cfrac{z}{q^2+\cfrac{z}{q^2+\cfrac{z}{q^3+...}}}}}} = \frac{\sum_{k=0}^{\infty}\frac{(-1)^k \cdot z^k \cdot q^{\frac{-k(k+1)}{2}}}{(q)_k}}{\sum_{k=0}^{\infty}\frac{(-1)^k \cdot z^k \cdot q^{\frac{-k(k-1)}{2}}}{(q)_k}}$$

where $|q|>1$. This is the Rogers Ramanujan continued fraction in another form.